\documentclass[12pt]{article}
\usepackage{latexsym,amssymb,amsmath}
\begin{document}

\baselineskip=18pt

\newcommand{\rd}{\mbox{Rad}}
\newcommand{\kn}{\mbox{ker}}
\newcommand{\psp}{\vspace{0.4cm}}
\newcommand{\pse}{\vspace{0.2cm}}
\newcommand{\ptl}{\partial}
\newcommand{\dlt}{\delta}
\newcommand{\Dlt}{\Delta}
\newcommand{\si}{\sigma}
\newcommand{\cta}{\theta}
\newcommand{\al}{\alpha}
\newcommand{\be}{\beta}
\newcommand{\G}{\Gamma}
\newcommand{\g}{\gamma}
\newcommand{\lmd}{\lambda}
\newcommand{\td}{\tilde}
\newcommand{\vf}{\varphi}
\newcommand{\ad}{\mbox{ad}}
\newcommand{\stl}{\stackrel}
\newcommand{\ol}{\overline}
\newcommand{\es}{\epsilon}
\newcommand{\vsi}{\varsigma}
\newcommand{\ves}{\varepsilon}
\newcommand{\la}{\langle}
\newcommand{\ra}{\rangle}
\newcommand{\vt}{\vartheta}
\newcommand{\wt}{\mbox{wt}\:}
\newcommand{\sym}{\mbox{sym}}
\newcommand{\for}{\mbox{for}}
\newcommand{\rta}{\rightarrow}
\newcommand{\der}{\mbox{Der}}
\newcommand{\mbb}{\mathbb}

\begin{center}{\Large \bf Central Simple Poisson Algebras}\end{center}

\begin{center}{\large Yucai Su$^{\ast}$ and  Xiaoping Xu$^{\dag}$}\end{center}

* Department of Applied Mathematics, Shanghai Jiaotong University, 1954 Huashan Road, Shanghai 200030, P. R. China.

\dag Department of Mathematics, The Hong Kong University of Science \& Technology, Clear Water Bay, Kowloon, Hong Kong, P. R. China\footnote{Research supported by Hong Kong RGC Competitive Earmarked Research Grant HKUST6133/00p}

\vspace{0.3cm}

\begin{abstract}{In this paper, we determine the isomorphism classes of the central simple Poisson algebras introduced earlier by the second author. The Lie algebra structures of these Poisson algebras are in general not finitely-graded.}\end{abstract}

\section{Introduction}

A {\it Poisson algebra} is a vector space ${\cal A}$ with two algebraic operations $\cdot$ and $[\cdot,\cdot]$ such that $({\cal A},\cdot)$ forms a commutative associative algebra, $({\cal A},[\cdot,\cdot])$ forms a Lie algebra and the following compatibility condition holds:
$$[u, v\cdot w]=[u,v]\cdot w+v\cdot[u,w]\qquad\for\;\;u,v,w\in {\cal A}.\eqno(1.1)$$
Poisson algebras are fundamental algebraic structures on phase spaces in classical mechanics. They are also the main objects in symplectic geometry. However, the structure theory of Poisson algebras does not seem to be well developed. 

 Let $({\cal A},\cdot, [\cdot,\cdot])$ be a Poisson algebra. Define
$$\mbox{Center}\:{\cal A}=\{u\in{\cal A}\mid [u,v]=0\;\for\;v\in{\cal A}\}.\eqno(1.2)$$
Then $\mbox{Center}\:{\cal A}$ is a Lie ideal of $({\cal A},[\cdot,\cdot])$. Form the quotient Lie algebra
$${\cal H}={\cal A}/(\mbox{Center}\:{\cal A})\eqno(1.3)$$
of $({\cal A},[\cdot,\cdot])$. We call the Poisson algebra ${\cal A}$ {\it central simple} if  $[{\cal H},{\cal H}]$ is a simple Lie algebra. A Lie algebra ${\cal G}$ is called {\it finitely-graded} if ${\cal G}=\bigoplus_{\al\in\G}{\cal G}_{\al}$ is a $\G$-graded vector space for some abelian group $\G$ such that
$$\dim\:{\cal G}_{\al}<\infty,\;\;[{\cal A}_{\al},{\cal A}_{\be}]\subset{\cal A}_{\al+\be}\qquad\for\;\;\al,\be\in \G.\eqno(1.4)$$

Central simple Poisson algebras, whose Lie algebra structures are finitely-graded, have been studied by Kac [K1], [K2], Osborn [O], and Osborn and Zhao [OZ].
The second author of this paper [X] constructed a family of central simple Poisson algebras, whose Lie algebra structures are in general not finitely-graded. The aim of this paper is to determine the isomorphism classes of central simple Poisson algebras given in [X].

Throughout this paper, we denote by $\mbb{F}$ a field with characteristic 0. All the vector spaces (algebras) are assumed over $\mbb{F}$. Moreover, we denote
by $\mbb{Z}$ the ring of integers and by $\mbb{N}$ the additive semi-group of nonnegative integers. When the context is clear, we shall omit the symbol for associative algebraic operation in a product. For $m,n\in\mbb{N}$, we shall use the following notation of indices
$$\ol{m,n}=\left\{\begin{array}{ll}\{m,m+1,...,n\}&\mbox{if}\;m\le n
\\ \emptyset&\mbox{if}\;m>n.\end{array}\right.
\eqno(1.5)$$

A classical central simple Poisson algebra is a polynomial algebra $\mbb{F}[t_1,t_2,...,t_{2\ell}]$ in $2\ell$ variables with the Lie bracket
$$[f,g]=\sum_{i=1}^{\ell}(\ptl_{t_i}(f)\ptl_{t_{\ell+i}}(g)-\ptl_{t_{\ell+i}}(f)\ptl_{t_i}(g))\qquad\for\;\;f,g\in \mbb{F}[t_1,t_2,...,t_{2\ell}].\eqno(1.6)$$
Define the grading
$$ (\mbb{F}[t_1,t_2,...,t_{2\ell}])_n=\mbox{Span}\:\{t_1^{n_1}t_2^{n_2}\cdots t^{n_{2\ell}}_{2\ell}\mid n_i\in\mbb{N},\;\sum_{i=1}^{2\ell}n_i=n+2\}\eqno(1.7)$$
for $-2\leq n\in\mbb{Z}$. Then  $\mbb{F}[t_1,t_2,...,t_{2\ell}]$ is a finitely-graded Lie algebra with respect to the bracket in (1.6).

The central simple Poisson algebras constructed in [X] are as follows. Let ${\cal A}=\bigoplus_{\al\in \G}{\cal A}_{\al}$ be a $\G$-graded commutative associative algebra with an identity element for some torsion-free abelian group $\G$, namely,
$${\cal A}_{\al}{\cal A}_{\be}\subset {\cal A}_{\al+\be}\qquad\for\;\;\al,\be\in\G.\eqno(1.8)$$
Take a skew-symmetric $\mbb{Z}$-bilinear form $\phi: \G\times\G \rta \mbb{F}$ and a set $\{\ptl_1,...,\ptl_{2\ell}\}$  of mutually commutative grading-preserving derivations of ${\cal A}$, that is,
$$\ptl_i({\cal A}_{\al})\subset {\cal A}_{\al},\;\;\ptl_i(uv)=\ptl_i(u)v+u\ptl_i(v)\qquad\for\;\;i\in\ol{1,2\ell},\;\al\in\G,\;u,v\in{\cal A}.\eqno(1.9)$$
Pick any elements
$$\{\xi_1,\xi_2,...,\xi_{\ell}\}\subset {\cal A}_0\eqno(1.10)$$
such that
$$\ptl_i(\xi_j)=\ptl_{\ell+i}(\xi_j)=0\qquad\for\;\;i,j\in\ol{1,\ell},\;i\neq j.\eqno(1.11)$$
We define an algebraic operation $[\cdot,\cdot]$ on ${\cal A}$ by
$$[u,v]=\phi(\al,\be)uv+\sum_{i=1}^{\ell}\xi_i(\ptl_i(u)\ptl_{\ell+i}(v)-\ptl_{\ell+i}(u)\ptl_i(v))\qquad\for\;\;u\in{\cal A}_{\al},\;v\in{\cal A}_{\be}.\eqno(1.12)$$
Then $({\cal A},\cdot, [\cdot,\cdot])$ forms a Poisson algebra. A linear transformation $T$ on a vector space $V$ is called {\it locally-finite} if
$$\dim(\mbox{Span}\:\{T^n(u)\mid n\in\mbb{N}\})<\infty\qquad\for\;\;u\in V.\eqno(1.13)$$

The central simple Poisson algebras constructed in [X] are of the above form $({\cal A},\cdot,[\cdot,\cdot])$, where ${\cal A}$ is a certain semi-group algebra, $\xi_i$ for $i\in\ol{1,\ell}$ are invertible, $\phi$ is nondegenerate, 
$\ptl_j$ for $j\in\ol{1,2\ell}$ are locally-finite, and some other distinguishable conditions among $\{\phi,\ptl_j\mid j\in\ol{1,2\ell}\}$ hold. We refer to [SXZ] for the classification of derivation-simple algebras when the derivations are locally-finite. The Lie algebra $({\cal A},[\cdot,\cdot])$ is in general not finitely-graded. The central simple Poisson algebras studied in [K1], [O] and [OZ] are special cases of those in [X], where the Lie algebra structures are finitely-graded. An isomorphism $\cta$ from a Poisson algebra ${\cal A}$ to another Poisson algebra ${\cal A}'$ is a linear isomorphism such that
$$\cta(uv)=\cta(u)\cta(v),\;\;\cta([u,v])=[\cta(u),\cta(v)]\qquad\for\;\;u,v\in{\cal A}.\eqno(1.14)$$

In Section 2, we shall rewrite the presentations of central simple Poisson algebras given in [X] up to certain relatively obvious isomorphisms, which we call {\it normalized forms}. The isomorphism classes of the normalized  central simple Poisson algebras will be determined in Section 3.

\section{Normalized Forms}

In this section, we shall give more concrete presentations of the central simple Poisson algebras constructed in [X] up to certain relatively obvious isomorphisms.

Let
$$\vec\ell=(\ell_1,...,\ell_7)\in \mbb{N}^{\:7}.
\eqno(2.1)$$
Set
$$\iota_i=\ell_1+\ell_2+...+\ell_i\qquad \for\ \ i\in\ol{1,7}
\eqno(2.2)$$
and
 define index sets
$$ I_{i,j}=\ol{\iota_{i-1}+1,\iota_j}\qquad\for\;\;i,j\in\ol{1,7},\;i\leq j,
\eqno(2.3)$$
where we treat $\iota_{i-1}$ as zero if $i=1$. For later convenience, we denote
$$I_i=I_{i,i},\;\;I=I_{1,7},\;\;J=\ol{1,2\iota_7}.\eqno(2.4)$$
Define the map $\ol{\ }: J\rta J$ by
$$\ol p=\left\{\begin{array}{ll}p+\iota_7&\mbox{if \ }p\in\ol{1,\iota_7},\\ p-\iota_7&\mbox{if \ }p\in\ol{\iota_7+1,2\iota_7}.\end{array}\right.
\eqno(2.5)$$
Moreover, for any subset $K$ of $\ol{1,2\iota_7}$, we denote
$$\ol K=\{\ol p\,|\,p\in K\}.\eqno(2.6)$$
In particular, we have
$$J=I\cup\ol I.\eqno(2.7)$$
Furthermore, we set
$$
J_i=I_i\cup\ol I_i,\ \ J_{i,j}=I_{i,j}\cup\ol I_{i,j}.
\eqno(2.8)$$

We shall always write an element $\al$ of $\mbb{F}^{\:2\iota_7}$ in the form
$$
\al=(\al_1,\al_{_{\ol1}},...,\al_{\iota_7},\al_{_{\ol\iota_7}})\qquad\mbox{with}\;\;\al_p\in\mbb{F}.\eqno(2.9)$$
Moreover, for $\al\in\mbb{F}^{\:2\iota_7}$ and $K\subset J$, we denote by 
$$
\al_{_ K}=(\be_1,\be_{_{\ol1}},...,\be_{\iota_7},\be_{_{\ol\iota_7}})\ \mbox{ with }\
\be_p=0\mbox{ \ if \ }p\notin K\mbox{ \ and \ }\be_p=\al_p\mbox{ \ if \ }p\in K.
\eqno(2.10)$$
Furthermore, we set
$$a_{_{[p]}}=(0,...,0,\stl{p}{a},0,...,0)\qquad\for\;\;a\in\mbb{F},\;\;p\in J.\eqno(2.11)$$
When the context is clear, we also use $\al_K$ to denote the vector obtained from $\al\in\mbb{F}^{\:2\iota_7}$ by deleting all  the coordinate $\al_p$ with $p\in J\setminus K$. For instance,
$$\al_{_{\{1,3\}}}=(\al_1,\al_3)\in\mbb{F}^{\:2},\;\;\al_{_{\{1,\bar{2},\bar{3}\}}}=(\al_1,\al_{\bar{2}},\al_{\bar{3}})\in\mbb{F}^{\:3}.\eqno(2.12)$$

Take
$$\si_p=\left\{\begin{array}{ll}1_{[p]}\mbox{ or }1_{[p]}+1_{[\ol p]}&\mbox{if \ }p\in J_1,\\ 1_{[p]}+1_{[\ol p]}&\mbox{if \ }p\in J_{2,3},\\ 0 &\mbox{if \ }p\in J_{4,7}.\end{array}\right.\eqno(2.13)$$
Let $\G_1$ be an additive subgroup of $\mbb{F}^{\:2\iota_7}$ such that
$$\al_{_{I_7\cup\ol I_{5,7}}}=0\qquad\for\;\;\al\in \G_1\eqno(2.14)$$
and
$$\{\si_p,\,1_{[q]}\,|\,p\in I,q\in I_{4,6}\cup\ol I_4\}\subset\G_1,\ \
\mbb{F}1_{[p]}\cap\G_1\ne\{0\}\ \ \for\ \ p\in J_{1,3}.
\eqno(2.15)$$
Let $\G_0$ be a torsion-free abelian group. Define
$$
\G=\G_0\oplus \G_1,\eqno(2.16)$$
a direct sum of abelian groups, and denote an element in $\G$ by
$$
\vec\al=\al_0+\al\qquad\mbox{with}\;\; \al_0\in\G_0,\ \al\in\G_1.
\eqno(2.17)$$

Note that  $\mbb{N}^{\:2\iota_7}$ is an additive semi-subgroup of $\mbb{F}^{\:2\iota_7}$. We take
$${\cal J}=\{\vec i=(i_1,i_{\bar{1}},...,i_{\iota_7}, i_{_{\ol \iota_7}})\in\mbb{N}^{\:2\iota_7}\mid \vec i_{I_2\cup I_5\cup \ol I_{1,2}}=0\}.
\eqno(2.18)$$
Let ${\cal A}$ be a vector space with a basis
$\{x^{\vec\al,\vec i}\mid(\vec\al,\vec i)\in\G\times {\cal J}\}$, and define the multiplication
$$
x^{\vec\al,\vec i}\cdot x^{\vec\be,\vec j}=x^{\vec\al+\vec\be,\vec i+\vec j}
\qquad\for\ \ (\vec\al,\vec i),(\vec\be,\vec j)\in\G\times {\cal J}.\eqno(2.19)$$
Then ${\cal A}$ forms a commutative associative algebra with $1=x^{0,0}$
as the identity element. For convenience, we denote
$$x^{\vec\al}=x^{\vec\al,0},\ \ t^{\vec i}=x^{0,\vec i},\ \
t_p=t^{1_{[p]}}.\eqno(2.20)$$
In particular,
$$t^{\vec i}=\prod_{p\in J}t_p^{i_p}\;\;\mbox{and}\;\;x^{\vec\al,\vec i}=x^{\vec\al}t^{\vec i}.\eqno(2.21)$$

Define the derivations $\{\ptl_p,\ptl^*_p,\ptl_{t_p}\,|\,p\in J\}$ of ${\cal A}$  by
$$
\ptl_p=\ptl^*_p+\ptl_{t_p}\mbox{ \ and \ }
\ptl^*_p(x^{\vec\al,\vec i})=\al_p x^{\vec\al,\vec i},\ \
\ptl_{t_p}(x^{\vec\al,\vec i})=i_p x^{\vec\al,\vec i-1_{[p]}},
\eqno(2.22)$$
for $p\in J,\;(\vec\al,\vec i)\in\G\times{\cal J},$
where we adopt the convention that if a notion is not defined but technically
appears in an expression, we always treat it as zero; for instance,
$x^{\vec\al,-1_{[1]}}=0$ for any $\vec\al\in\G$.
In particular,
$$
\ptl^*_p=0,\ \ \ptl_{t_q}=0\qquad\for\ \ p\in
I_7\cup\ol I_{5,7},\ q\in I_2\cup I_5\cup \ol I_{1,2},\eqno(2.23)$$
by (2.14) and (2.18). We call the nonzero derivations $\ptl^*_p$ {\it grading operators}, the nonzero derivations $\ptl_{t_q}$ {\it down-grading operators}, and the derivations $\ptl_r^{\ast}+\ptl_{t_r}$ {\it mixed operators} if both $\ptl_r^{\ast}$ and $\ptl_{t_r}$ are not zero. The types of derivation pairs in the order of the groups   $\{(\ptl_p,\ptl_{\bar{p}})\mid p\in I_i\}$ for $i\in\ol{1,7}$ are
$$ (m,g),\;(g,g),\;(m,m),\;(m,m),\;(g,d),\;(m,d),\;(d,d),\eqno(2.24)$$
where ``m'' stands for mixed operators, ``g'' stands for grading operators and ``d'' stands for down-grading operators.

Let $\phi(\cdot,\cdot):\G\times\G\rta\mbb{F}$ be a skew-symmetric $\mbb{Z}$-bilinear
form such that
$$\si_p\in\rd _\phi\qquad \for\ \ p\in I\eqno(2.25)$$
and
$$\mbb{F}1_{[p]}\cap\rd_\phi\ne\{0\}\qquad \for\ \ p\in J_{1,3}.
\eqno(2.26)$$
Set
$$
{\cal A}_{\vec\al}={\rm Span}\{x^{\vec\al,\vec i}\,|\,\vec i\in{\cal J}\}\ \ \for\ \ \vec\al\in\G.
\eqno(2.27)$$
We define the following Lie bracket $[\cdot,\cdot]$ on ${\cal A}$:
$$
[u,v]=\sum_{p\in I}x^{\si_p}
(\ptl_p(u)\ptl_{_{\ol p}}(v)-\ptl_{_{\ol p}}(u)\ptl_p(v))
+(\phi(\al,\be)-\sum_{p\in I_4}(\al_p\be_{_{\ol p}}-\al_{_{\ol p}}\be_p))uv,
\eqno(2.28)$$
for $u\in{\cal A}_{\vec\al},v\in{\cal A}_{\vec\be}$. Then $({\cal A},\cdot,[\cdot,\cdot])$
forms a Poisson algebra.

The above Poisson algebras $({\cal A},\cdot,[\cdot,\cdot])$ are the normalized forms
of Poisson algebras constructed in [X]. Moreover,
$$\mbox{center}\:{\cal A}=\mbb{F}1_{\cal A}.\eqno(2.29)$$
The following theorem was proved in [X].
\psp

{\bf Theorem 2.1}. {\it If
$$\{\al_0\in\G_0\mid\phi(\al_0,\vec\be)=0\;\mbox{\it for}\;\vec\be\in\G\;\mbox{\it
with}\;\be_{_{J_{1,3}}}=0\}=\{0\},\eqno(2.30)$$
then $({\cal A},\cdot,[\cdot,\cdot])$ is central simple}.
\psp

Thus in rest of this paper,  we shall assume that (2.30) holds. Moreover, we denote the above Poisson algebra by
$${\cal P}(\vec \ell,\G,{\cal J},\si,\phi)=({\cal A},\cdot,[\cdot,\cdot]),\eqno(2.31)$$
in order to emphasize its dependence on the constructional ingredients, where
$$\si=\sum_{i=1}^{\iota_3}\si_i.\eqno(2.32)$$

\section{Isomorphism Classes}

In this section, we shall determine the isomorphism classes of the Poisson algebras of the form ${\cal P}(\vec \ell,\G,{\cal J},\si,\phi)$. We assume that $\mbb{F}$ is an algebraically closed field.

Consider a specific Poisson algebra ${\cal P}(\vec \ell,\G,{\cal J},\si,\phi)$.
By re-indexing the index set $J_1$ if necessary, we may assume
$$\si_p=\left\{\begin{array}{ll}
1_{[p]}&\mbox{if \ }p\in\ol{1,\ell_0}\mbox{ or }\ol p\in\ol{1,\ell_0},\\
1_{[p]}+1_{[\ol p]}&\mbox{if \ }p\in\ol{\ell_0+1,\ell_1}\mbox{ or }\ol p
\in\ol{\ell_0+1,\ell_1},\end{array}\right.\eqno(3.1)$$
for some $\ell_0\in\ol{0,\ell_1}$.
By (2.14), (2.18), (2.22) and (2.23), we can rewrite (2.28) as follows.
\begin{eqnarray*}[x^{\vec\al,\vec i},x^{\vec\be,\vec j}]
&=&\sum_{p\in I_{1,3}}(\al_p\be_{_{\ol p}}-\al_{_{\ol p}}\be_p)
x^{\si_p+\vec\al+\vec\be,\vec i+\vec j}+\sum_{p\in I_{3,6}}(\al_pj_{_{\ol p}}-i_{_{\ol p}}\be_p)x^{\si_p+\vec\al+\vec\be,\vec i+\vec j-1_{[\ol p]}}\\ &&+
\sum_{p\in I_1\cup I_{3,4}}(i_p\be_{_{\ol p}}-j_p\al_{_{\ol p}})
x^{\si_p+\vec\al+\vec\be,\vec i+\vec j-1_{[p]}}+\phi(\vec\al,\vec\be)x^{\vec\al+\vec\be,\vec i+\vec j}\\ &&+\sum_{p\in I_{3,4}\cup I_{6,7}}(i_pj_{_{\ol p}}-i_{_{\ol p}}j_p)x^{\si_p+\vec\al+\vec\be,\vec i+\vec j-1_{[p]}-1_{[\ol p]}}\hspace{5.1cm}(3.2)\end{eqnarray*}
for $(\vec\al,\vec i),(\vec\be,\vec j)\in\G\times{\cal J}$.

Denote by $M_{m\times n}(\mbb{F})$ the space of $m\times n$ matrices with entries in $\mbb{F}$ and by $GL_m(\mbb{F})$ the group of $m\times m$ invertible matrices with entries in $\mbb{F}$. Define
$$G_p=\left\{ \left(\begin{array}{cc}1&0\\ a&b\end{array}\right)\mid a,b\in\mbb{F},\;b\neq 0\right\}\qquad \for\;\;p\in\ol{1,\ell_0},\eqno(3.3)$$
$$G_q=\left\{\left(\begin{array}{cc}a&0\\ 1-a&1\end{array}\right)\mid 0\neq a\in\mbb{F}\right\}\qquad \for\;\;q\in\ol{\ell_0+1,\ell_1},\eqno(3.4)$$
$$G_r=\left\{ \left(\begin{array}{cc}a+b&a\\ 1-a-b&1-a\end{array}\right)\mid a,b\in\mbb{F},\;b\neq 0\right\}\qquad \for\;\;r\in I_{2,3}.\eqno(3.5)$$
Set
$$S_m={\rm diag}\left(\left(\begin{array}{rr}0&-1\\ 1&0\end{array}\right),...,\left(\begin{array}{rr}0&-1\\ 1&0\end{array}\right)\right)\in GL_{2m}(\mbb{F})\eqno(3.6)$$
for $m\in\mbb{N}$. Let
$$SP_{2\ell_4}(\mbb{F})=\{ A\in GL_{2\ell_4}(\mbb{F})\mid A^TS_{\ell_4}A=S_{\ell_4}\},\eqno(3.7)$$
where the up-index ``T'' means the transpose of matrices. Denote by ${\bf 1}_m$ the $m\times m$ identity matrix.
Moreover, we define  $F$ to be the group of invertible matrices of the form
$$\left(\begin{array}{cccccc} A_1&0&C&0&0&0\\ B&A_2&D&0&0&0\\  0&0&A_3&0&0&0 \\ 0&0&0&{\bf 1}_{\ell_5}&0&0\\ -(A_3^T)^{-1}C^TS_{\ell_4}A_1& 0&-\frac{1}{2}(A_3^T)^{-1}C^TS_{\ell_4}C&0&(A_3^T)^{-1}&0\\ 0&0&0&0&0&{\bf 1}_{2\ell_7}\end{array}\right),\eqno(3.8)$$
where 
$$  A_1\in SP_{2\ell_4}(\mbb{F}),\;\;A_2\in GL_{\ell_5}(\mbb{F}),\;\;A_3\in GL_{\ell_6}(\mbb{F}),\eqno(3.9)$$
$$B\in M_{\ell_5\times 2\ell_4}(\mbb{F}),\;\;C\in M_{2\ell_4\times \ell_6}(\mbb{F}),\;\;D\in M_{\ell_5\times \ell_6}(\mbb{F}).\eqno(3.10)$$
Now we define the group
$$G'=\{\mbox{diag}(g_1,...,g_{\iota_3},f)\in GL_{2\iota_7}(\mbb{F})\mid g_i\in G_i,\;f\in F\}.\eqno(3.11)$$

Let ${\cal S}_I$ be the permutation group on the index set $I$ (cf. (2.4)). Define the subgroup
\begin{eqnarray*}\hspace{2cm}{\cal S}&=&\{\nu\in {\cal S}_I\mid \nu(\ol{1,\ell_0})=\ol{1,\ell_0},\;\nu(\ol{\ell_0+1,\ell_1})=\ol{\ell_0+1,\ell},\\ & &\nu(I_2)=I_2,\;\nu(I_3)=I_3,\;\nu|_{I_{4,7}}=\mbox{Id}_{I_{4,7}}\}.\hspace{4.8cm}(3.12)\end{eqnarray*}
Moreover, for $g\in G'$ and $\nu\in{\cal S}$, we define the group automorphism $g_{\nu}$ of $\mbb{F}^{\:2\iota_7}$ by
\begin{eqnarray*}& &((g_{\nu}(\al))_{\nu(1)},(g_{\nu}(\al))_{_{\ol{\nu(1)}}},...,(g_{\nu}(\al))_{\nu(\iota_4)},(g_{\nu}(\al))_{_{\ol{\nu(\iota_4)}}},(g_{\nu}(\al))_{_{I_{5,6}}},(g_{\nu}(\al))_{_{\ol I_{5,6}\cup J_7}})\\&=&(\al_{_{J_{1,4}}},\al_{_{I_{5,6}}},\al_{_{\ol I_{5,6}\cup J_7}})g\hspace{10.1cm}(3.13)\end{eqnarray*}
(cf. (2.12)), where $\al\in \mbb{F}^{\:2\iota_7}$ and the multiplication in the above is the vector-matrix multiplication. Define
$$G=\{g_{\nu}\mid g\in G',\;\nu\in{\cal S}\}.\eqno(3.14)$$
Then $G$ is a subgroup of additive automorphims of $\mbb{F}^{\:2\iota_7}$. 
  
Let ${\cal P}(\vec \ell',\G',{\cal J}',\si',\phi')$ be another Poisson algebra defined in last section. We shall add a prime on all the  constructional ingredients related to ${\cal P}(\vec \ell',\G',{\cal J}',\si',\phi')$; for instance, ${\cal A}',\;\ell_i',\;\iota_i',$ etc.
\psp

{\bf Theorem 3.1}. {\it The Poisson algebra ${\cal P}(\vec \ell,\G,{\cal J},\si,\phi)$ is isomorphic to the Poisson algebra ${\cal P}(\vec \ell',\G',{\cal J}',\si',\phi')$ if and only if $(\ell_0,\vec\ell)=(\ell'_0,\vec\ell')$ and there exists a group isomorphism $\tau$ from $\G$ to $\G'$ of the form
$$\tau(\al_0+\al)=\tau_0(\al_0)+\tau_1(\al)+\tau_2(\al)\qquad\for\;\;\al_0\in\G_0,\;\al\in \G_1,\eqno(3.15)$$
where $\tau_0:\G_0\rta\G'_0$ is a group isomorphism, $\tau_1:\G_1\rta\G'_0$ is a group homomorphism and $\tau_2\in G$ such that $\tau_2(\G_1)=\G_1'$ and 
$$\phi'(\tau_0(\al_0),\tau_0(\be_0))=\phi(\al_0,\be_0),\;\;\phi'(\tau_0(\al_0),\tau_2(\al))+\phi'(\tau_0(\al_0),\tau_1(\al))=\phi(\al_0,\al),\eqno(3.16)$$
$$\phi'(\tau_1(\al),\tau_1(\be))+\phi'(\tau_1(\al),\tau_2(\be))
+\phi'(\tau_2(\al),\tau_1(\be))+\phi'(\tau_2(\al),\tau_2(\be))=\phi(\al,\be),
\eqno(3.17)$$
for  $\al_0,\be_0\in\G_0$ and $\al,\be\in\G_1$.}
\psp

{\it Proof}. Recall that ${\cal A}$ is a commutative associative algebra (cf. (2.19)). By  (2.21) and  (2.28), we have
\begin{eqnarray*}\hspace{2cm}[x^{\vec\al,\vec i},x^{\vec\be,\vec j}]&=&
[x^{\vec\al},x^{\vec\be}]t^{\vec i+\vec j}+\sum_{p\in J}(i_p x^{\vec\al}[t_p,x^{\vec\be}]-j_p[t_p,x^{\vec\al}]x^{\vec \be})t^{\vec i+\vec j-1_{[p]}}\\&&
+\sum_{p,q\in J}i_pj_qx^{\vec\al+\vec\be}[t_p,t_q]t^{\vec i+\vec j-1_{[p]}-1_{[q]}}\hspace{4.4cm}(3.18)\end{eqnarray*}
for  $\vec\al,\vec\be\in\G$ and $\vec i,\vec j\in{\cal J}$.
Thus an associative algebra isomorphism $\cta:{\cal A}\rta{\cal A}'$
is a Poisson algebra isomorphism if and only if it satisfies
$$\cta([x^{\vec\al},x^{\vec\be}])=[\cta(x^{\vec\al}),\cta(x^{\vec\be})],\;\;
\cta([t_p,x^{\vec\be}])=[\cta(t_p),\cta(x^{\vec\be})],\;\;\cta([t_p,t_q])=[\cta(t_p),\cta(t_q)]\eqno(3.19)$$
for all $\vec\al,\vec\be\in\G$ and  $p,q\in J\setminus(I_2\cup I_5\cup\ol I_{1,2})$ (cf. (2.18)).
For convenience, we denote
$$\ol t=(-t_{\ol1},t_1,-t_{\ol2},t_2,...,-t_{\ol\iota_7},t_{\iota_7}).
\eqno(3.20)$$
For a subset $K$ of $J$, we denote by  $\ol t_{_K}$ the vector obtained from $\ol t$ by deleting $-t_{\ol p}, t_q$ for $\ol p,q\in J\setminus K$. For instance,
$$\ol t_{\{\ol 1,\ol 2,2\}}=(-t_{\ol 1},-t_{\ol 2},t_2).\eqno(3.21)$$
\psp

``$\Longleftarrow$'' First we prove the sufficiency. Assume that $(\ell_0,\vec\ell)=(\ell'_0,\vec\ell')$ and there exist a group isomorphism $\tau$ from $\G$ to $\G'$ of the form (3.15) such that (3.16) and (3.17) hold.
\pse

By (3.14), 
$$\tau_2=g_{\nu}\;\;\mbox{with}\;\;g=\mbox{diag}(g_1,...,g_{\iota_3},f)\in G',\;\nu\in {\cal S}_I,\eqno(3.22)$$
where we write
$$g_p=\left(\begin{array}{cc}1&0\\ a_p&b_p\end{array}\right),\;\;g_q=\left(\begin{array}{cc}b_q&0\\ 1-b_q&1\end{array}\right),\;\;g_r=\left(\begin{array}{cc}a_r+b_r&a_r\\ 1-a_r-b_r&1-a_r\end{array}\right)\eqno(3.23)$$
for $p\in\ol{1,\ell_0},\;q\in\ol{\ell_0+1,\ell_1}$ and $r\in I_{2,3}$. 
In ${\cal A}'$, we define 
$$(s_p,-s_{_{\ol p}})=b_p(t'_{\nu(p)},-t'_{_{\ol{\nu(p)}}})g_p^{-1}\qquad\for\ \ p\in I_{1,3},\eqno(3.24)$$
$$(\ol s_{_{J_4\cup \ol{I}_{5,6}}},\ol s_{_{I_{5,6}\cup J_7}})=(\ol t'_{_{J_4\cup \ol{I}_{5,6}}},\ol t'_{_{I_{5,6}\cup J_7}})f^{-1}\eqno(3.25)$$
(cf. (2.8), (2.18), (3.8), (3.20) and (3.21)).

Let
$$\Dlt_0=\sum_{p=1}^{\iota_3}\mbb{Z}\si_p\eqno(3.26)$$
be the subgroup of $\G$ generated by $\{\si_i\mid i\in\ol{1,\iota_3}\}$ (cf. (2.13)). Define $\chi:\Dlt_0\rta\mbb{F}^{\times}$ to be the homomorphism from additive group to the multiplication group of nonzero elements of $\mbb{F}$ determined by
$$\chi(\si_p)=b_p\ \ \for\ \ p\in I_{1,3}.\eqno(3.27)$$
We want to prove that  $\chi$ can be extended to a homomorphism $\chi: \G\rta \mbb{F}^{\times}$. Suppose that $\Dlt$ is a maximal subgroup of $\G$ containing $\Dlt_0$ such that $\chi$ can be extended to a homomorphism $\chi: \Dlt\rta \mbb{F}^{\times}$. Assume $\Dlt\neq \G$. We take an element $\vec \al\in\G\setminus \Dlt$. Set
$$\Dlt'=\mbb{Z}\vec\al+\Dlt.\eqno(3.28)$$
If $\mbb{Z}\vec \al\bigcap\Dlt=\{0\}$, then we extend $\chi$ by
$$\chi(m\vec\al+\vec\be)=\chi(\vec\be)\qquad\for\;\;m\in\mbb{Z},\;\be\in\Dlt.\eqno(3.29)$$
If  $\mbb{Z}\vec \al\bigcap\Dlt=\mbb{Z}n\vec \al$, we take an $n$th root $a$ of
$\chi(n\vec\al)$ (recall that $\mbb{F}$ is algebraically closed) and  extend $\chi$ by
$$\chi(m\vec\al+\vec\be)=a^m\chi(\vec\be)\qquad\for\;\;m\in\mbb{Z},\;\be\in\Dlt.\eqno(3.30)$$
It is straightforward to verify that $\chi: \Dlt'\rta \mbb{F}^{\times}$ is a group homomorphism. This leads a contradiction to the maximality of $\Dlt$. So 
$\chi$ can be extended to a homomorphism $\chi: \G\rta \mbb{F}^{\times}$. Take any such extension.

Recall that  we add prime on the constructional ingredients related to ${\cal P}(\vec \ell',\G',{\cal J}',\si',\phi')$. Now we define the associative algebra  isomorphism $\cta: {\cal A}\rta {\cal A}'$ by
$$\cta(x^{\vec \al,\vec i})=\chi(\vec\al){x'}^{\tau(\vec\al)}s^{\vec i}\qquad \for\;\;(\vec \al,\vec i)\in \G\times{\cal J}\eqno(3.31)$$
(cf. (2.21), (3.24) and (3.25)). Moreover, (3.16) and (3.17) guarantee
$$\phi'(\tau(\vec\al),\tau(\vec\be))=\phi(\vec\al,\vec\be)\qquad \for\ \ \vec\al,\vec\be\in\G.\eqno(3.32)$$
Assume that $\tau(\si_p)=\vec\be'=(\be_0',\be')$  for $p\in I_{1,3}$.  By (3.3)-(3.13), we obtain that $\be'=\tau_2(\si_p)=\si'_{\nu(p)}$.
Since $\si_p\in\rd_\phi,\si'_{\nu(p)}\in\rd_{\phi'}$
(cf. (2.25)), $\be_0'\in\rd_{\phi'}$ by (3.32). Moreover, by our assumption of (2.30),  $\be_0'=0$. Hence
$$\tau(\si_p)=\si'_{\nu(p)}\qquad \for\ \ p\in I_{1,3}.\eqno(3.33)$$

By (2.10) and (3.2), we have
$$[x^{\vec\al},x^{\vec\be}]=\sum_{p\in I_{1,3}}
\left|\begin{array}{c}\al_{\{p,\ol p\}}\\\be_{\{p,\ol p\}}\end{array}\right|x^{\si_p+\vec\al+\vec\be}+\phi(\vec\al,\vec\be)x^{\vec\al+\vec\be},\eqno(3.34)$$
where
$$\left|\begin{array}{c}\al_{\{p,\ol p\}}\\\be_{\{p,\ol p\}}\end{array}\right|
=\left|\begin{array}{cc}\al_p&\al_{_{\ol p}}\\\be_p&\be_{_{\ol p}}\end{array}\right|\eqno(3.35)$$
is a $2\times2$ determinant. Moreover, by (3.2) and (3.31)-(3.33), we get
\begin{eqnarray*}[\cta(x^{\vec\al}),\cta(x^{\vec\be})]&=&\chi(\vec\al)\chi(\vec\be)(\sum_{p\in I_{1,3}}\left|\begin{array}{c}(\tau_2(\al))_{\{p,\bar{p}\}}\\ (\tau_2(\be))_{\{p,\ol p\}}\end{array}\right|
{x'}^{\si'_{\nu(p)}+\tau(\vec\al)+\tau(\vec\be)}\\ & &+
\phi'(\tau(\vec\al),\tau(\vec\be)){x'}^{\tau(\vec\al)+\tau(\vec\be)})\\ &=&\cta(
\sum_{p\in I_{1,3}}
\left|\begin{array}{c}\al_{\{p,\ol p\}}\\\be_{\{p,\ol p\}}\end{array}\right|x^{\si_p+\vec\al+\vec\be}+\phi(\vec\al,\vec\be)x^{\vec\al+\vec\be})\hspace{4.4cm}(3.36)\end{eqnarray*}
 because (3.23) implies the determinant of $g_p$ is $|g_p|=b_p=\chi(\si_p)$
(cf. (3.27)) and 
\begin{eqnarray*}\hspace{2cm}\chi(\vec\al)\chi(\vec\be)\left|\begin{array}{c}(\tau_2(\al))_{\{p,\bar{p}\}}\\ (\tau_2(\be))_{\{p,\ol p\}}\end{array}\right|&
=&\chi(\vec\al+\vec\be)\left|\begin{array}{c}\al_{\{p,\ol p\}}\\\be_{\{p,\ol p\}}\end{array}\right|\cdot|g_p|\\&=&\chi(\si_p+\vec\al+\vec\be)\left|\begin{array}{c}\al_{\{p,\ol p\}}\\\be_{\{p,\ol p\}}\end{array}\right|.\hspace{3.1cm}(3.37)\end{eqnarray*}

Recall the notation (2.10) and ${\cal J}\subset\mbb{F}^{\:2\iota_7}$ (cf. (2.18)). As $1\times2$ matrices,
$$[\ol t_{\{p,\ol p\}},x^{\vec\al}]=([-t_{_{\ol p}},x^{\vec\al}],[t_p,x^{\vec\al}])=\al_{\{p,\ol p\}}x^{\si_p+\vec\al}\qquad\for\ \ p\in I\eqno(3.38)$$
(cf. (2.4), (3.20) and (3.21)). Now we verify the second equation in (3.19). If $p\in I_1$, we have 
$$s_p=\left\{\begin{array}{ll}t'_{\nu(p)}&\mbox{if}\;p\in\ol{1,\ell_0},\\ b_pt'_{\nu(p)}&\mbox{if}\;p\in\ol{\ell_0+1,\ell_1}\end{array}\right.\eqno(3.39)$$
by (3.24). Moreover, by (3.13), (3.22) and (3.23),
$$(\tau_2(\al))_{_{\ol{\nu(p)}}}=\left\{\begin{array}{ll}b_p\al_{\bar{p}}&\mbox{if}\;p\in\ol{1,\ell_0},\\ \al_{\bar{p}}&\mbox{if}\;p\in\ol{\ell_0+1,\ell_1}.\end{array}\right.\eqno(3.40)$$
Thus
\begin{eqnarray*}\hspace{2cm}[\cta(t_p),\cta(x^{\vec\al})]&=&\chi(\vec\al)[s_p,{x'}^{\tau(\vec\al)}]\\ &
=&\chi(\vec\al)\left\{\begin{array}{ll}[t'_{\nu(p)},{x'}^{\tau(\vec\al)}]&\mbox{if}\;p\in\ol{1,\ell_0},\\ b_p[t'_{\nu(p)},{x'}^{\tau(\vec\al)}]&\mbox{if}\;p\in\ol{\ell_0+1,\ell_1}\end{array}\right.
\\ &=&\chi(\vec\al)\left\{\begin{array}{ll}(\tau_2(\al))_{_{\ol{\nu(p)}}}{x'}^{\tau(\vec\al)+\si'_{\nu(p)}}&\mbox{if}\;p\in\ol{1,\ell_0},\\ b_p(\tau_2(\al))_{_{\ol{\nu(p)}}}{x'}^{\tau(\vec\al)+\si'_{\nu(p)}}&\mbox{if}\;p\in\ol{\ell_0+1,\ell_1}\end{array}\right.\\&=&\chi(\vec\al)b_p\al_{\bar{p}}{x'}^{\tau(\vec\al)+\si'_{\nu(p)}}\\ &=&\chi(\vec\al+\si_p)\al_{\bar{p}}{x'}^{\tau(\vec\al)+\si'_{\nu(p)}}\\ &=&\cta(\al_{\bar{p}}x^{\vec\al+\si_p})\\ &=&\cta([t_p,x^{\vec \al}]).\hspace{7.8cm}(3.41)\end{eqnarray*}
For $p\in I_3$, by (3.13) and (3.20)-(3.23), we have
\begin{eqnarray*}\hspace{3cm}[\cta(\ol t_{\{p,\ol p\}}),\cta(x^{\vec\al})]&=&\chi(\vec\al)[\ol s_{\{p,\ol p\}},{x'}^{\tau(\vec\al)}]\\ &=&\chi(\vec\al)[\ol t'_{\{\nu(p),\ol {\nu(p)}\}},{x'}^{\tau(\vec\al)}]g_p^{-1}\\ &=&\chi(\vec\al)b_p(\tau_2(\al))_{\{\nu(p),\ol{\nu(p)}\}}g^{-1}_p{x'}^{\si'_{\nu(p)}+\tau(\vec\al)}\\ &=&\chi(\vec\al)b_p\al_{\{p,\ol p\}} {x'}^{\si'_{\nu(p)}+\tau(\vec\al)}\\ &=& \chi(\vec\al+\si_p)\al_{\{p,\ol p\}} {x'}^{\si'_{\nu(p)}+\tau(\vec\al)}\\ &=&\cta(\al_{\{p,\ol p\}} x^{\vec\al+\si_p})\\ &=&\cta([\ol t_{\{p,\bar{p}\}},x^{\vec\al}]).\hspace{5.8cm}(3.42)\end{eqnarray*}
Furthermore, by (3.25), we get
\begin{eqnarray*}[(\cta(\ol t_{_{J_4\cup \ol I_{5,6}}}), \cta(\ol t_{_{I_{5,6}\cup J_7}})),\cta(x^{\vec\al})]&=&
\chi(\vec\al)[(\ol s_{_{J_4\cup \ol I_{5,6}}},\ol s_{_{I_{5,6}\cup J_7}}),
{x'}^{\tau(\vec\al)}]\\ &=&\chi(\vec\al)[(\ol t'_{_{J_4\cup \ol I_{5,6}}},\ol t'_{_{I_{5,6}\cup J_7}}),{x'}^{\tau(\vec\al)}]f^{-1}\\ &=&\chi(\vec\al)((\tau_2(\al))_{_{J_4\cup \ol I_{5,6}}},(\tau_2(\al))_{_{I_{5,6}\cup J_7}})f^{-1}{x'}^{\tau(\vec\al)}\\ &=&\chi(\vec\al)(\al_{_{J_4\cup \ol I_{5,6}}},\al_{_{I_{5,6}\cup J_7}}) {x'}^{\tau(\vec\al)}\\ &=& \cta((\al_{_{J_4\cup \ol I_{5,6}}},\al_{_{I_{5,6}\cup J_7}}) x^{\vec\al})\\ &=&\cta([(\ol t_{_{J_4\cup \ol I_{5,6}}},\ol t_{_{I_{5,6}\cup J_7}}),x^{\vec\al}]).\hspace{4.5cm}(3.43)\end{eqnarray*}
This proves the second equation of (3.19).

Set
$$\Psi=\left(\begin{array}{cccccc} S_{\ell_4}&0&0&0&0&0\\ 0&0&0&0&0&0\\  0&0&0&0&{\bf 1}_{\ell_6}&0 \\ 0&0&0&0&0&0\\ 0&0&-{\bf 1}_{\ell_6}&0&0&0\\ 0&0&0&0&0&S_{\ell_7}\end{array}\right)\eqno(3.44)$$
(cf. (3.6)), where ${\bf 1}_{\ell_6}$ is the $\ell_6\times\ell_6$ identity matrix.
Then we have
\begin{eqnarray*}\hspace{2cm}& &[\cta(\ol t_{_{J_4\cup \ol{I}_{5,6}}},\ol t_{_{I_{5,6}\cup J_7}})^T,\cta(\ol t_{_{J_4\cup \ol{I}_{5,6}}},\ol t_{_{I_{5,6}\cup J_7}})]\\ &=&[(\ol s_{_{J_4\cup \ol{I}_{5,6}}},\ol s_{_{I_{5,6}\cup J_7}})^T,(\ol s_{_{J_4\cup \ol{I}_{5,6}}},\ol s_{_{I_{5,6}\cup J_7}})]\\ &=&[(f^{-1})^T(\ol t'_{_{J_4\cup \ol{I}_{5,6}}},\ol t'_{_{I_{5,6}\cup J_7}})^T,(\ol t'_{_{J_4\cup \ol{I}_{5,6}}},\ol t'_{_{I_{5,6}\cup J_7}})f^{-1}]\\ &=& (f^{-1})^T\Psi f^{-1}=\Psi\\& =&\cta([(\ol t_{_{J_4\cup \ol{I}_{5,6}}},\ol t_{_{I_{5,6}\cup J_7}})^T,(\ol t_{_{J_4\cup \ol{I}_{5,6}}},\ol t_{_{I_{5,6}\cup J_7}})])\hspace{5.2cm}(3.45)\end{eqnarray*}
by (3.2), (3.7), (3.8) and (3.25). Similarly, we can prove
$$[\cta(\ol t_{_{J_3}})^T,\cta(\ol t_{_{J_3}})]=\cta([\ol t_{_{J_3}}^T,\ol t_{_{J_3}}])\eqno(3.46)$$
by (3.23) and (3.24). The other identities hold trivially. Therefore, the last equation in (3.19) holds. This proves that the map $\cta: {\cal A}\rta {\cal A}'$ in (3.31) is a Poisson algebra isomorphism from ${\cal P}(\vec \ell,\G,{\cal J},\si,\phi)$ to ${\cal P}(\vec \ell',\G',{\cal J}',\si',\phi')$.
\psp

``$\Longrightarrow$'' We assume that there exits a Poisson algebra isomorphism $\cta: {\cal P}(\vec \ell,\G,{\cal J},\si,\phi)\rta {\cal P}(\vec \ell',\G',{\cal J}',\si',\phi')$.
\pse

First, we make the following conventions. If a subset of ${\cal A}$ is defined,
then we take the definition of the corresponding subset of ${\cal A}'$ for granted. If a property about ${\cal P}(\vec \ell,\G,{\cal J},\si,\phi)$ is given, the same property also holds for ${\cal P}(\vec \ell',\G',{\cal J}',\si',\phi')$, without description.

Consider the commutative associative algebra structure $({\cal A},\cdot)$. It can be proved by taking an order on $\G$ that
$$\{u\in{\cal A}\,|\,u\mbox{ is invertible in }({\cal A},\cdot)\}
=(\bigcup_{\vec\al\in\G}\mbb{F} x^{\vec\al})\setminus\{0\}.
\eqno(3.47)$$
The same statement holds for ${\cal A}'$. Hence there exists a bijection $\tau:\G\rta \G'$ such that
$$\cta(x^{\vec\al})=\chi(\vec\al) {x'}^{\tau(\vec\al)}\qquad\for\ \ \vec\al\in\G\mbox{ \ and some \ }\chi({\vec\al})\in\mbb{F}\setminus\{0\}.\eqno(3.48)$$
Using the fact $x^{\vec\al}\cdot x^{\vec\be}=x^{\vec\al+\vec\be}$, we obtain that
$\tau$ is an group isomorphism and $\chi$ is a homomorphism $\chi$ from $\G$ to $\mbb{F}^{\times}$.

Denote 
$${\cal A}^{(0)}={\rm Span}\{x^{\vec\al}\mid\vec\al\in\G\},\eqno(3.49)$$
 the group algebra $\G$. For $\vec\al\in\G$ and $\mu\in\mbox{Hom}_{\mbb{Z}}(\G, \mbb{F})$ the space of additive group homomorphisms,  
we define a derivation of ${\cal A}^{(0)}$ by
$$d_{\vec\al,\mu}(x^{\vec\be})=\mu(\vec\be)x^{\vec\al+\vec\be}\qquad\for\;\;\vec\be\in\G.\eqno(3.50)$$
Denote by $\der{\cal A}^{(0)}$ the derivation algebra
of ${\cal A}^{(0)}$ and
$$(\der{\cal A}^{(0)})_{\vec\al}={\rm Span}\{d_{\vec\al,\mu}\mid \mu\in\mbox{Hom}_{\mbb{Z}}(\G, \mbb{F})\}.\eqno(3.51)$$
It can be proved by taking an order on $\G$ that
$$
\der{\cal A}^{(0)}=\bigoplus_{\vec\al\in\G}(\der{\cal A}^{(0)})_{\vec\al}.\eqno(3.52)$$
Note that a nonzero derivation $d_{\vec\al,\mu}$ of ${\cal A}^{(0)}$ is
locally finite (cf. (1.13)) if and only if it is a homogeneous derivation of degree $\vec\al=0$.

We shall prove the necessity of the theorem by establishing seven claims (steps).  
\psp

{\bf Claim 1}.  We have
$$\{\vec\al\in\G\mid\ad_{x^{\vec\al}}|_{{\cal A}^{(0)}}=0\}=\{\vec\al\in\rd_\phi\,|\,\al_{J_{1,3}}=0\}\eqno(3.53)$$
(cf. (2.10)),
$$\{\vec\al\in\G\mid\ad_{x^{\vec\al}}|_{{\cal A}^{(0)}}\ne0\mbox{ is locally-finite}\}=\{-\si_p\,|\,p\in I_{1,3}\},\eqno(3.54)$$
$$\{p\in I_{1,3}\mid \ad_{x^{-\si_p}}|_{{\cal A}^{(0)}}\mbox{ is diagonalizable} \}=\ol{1,\ell_0}\cup I_2,\eqno(3.55)$$
$$\{p\in I_{1,3}\mid \ad_{x^{-\si_p}}|_{{\cal A}^{(0)}}\mbox{ is not diagonalizable}\}=\ol{\ell_0+1,\ell_1}\cup I_3.\eqno(3.56)$$
\pse

For a given $\al\in \G$, by (3.34), we have
$$\ad_{x^{\vec\al}}|_{{\cal A}^{(0)}}=\sum_{p\in I_{1,3}}d_{\si_p+\vec\al,\mu_p}+d_{\vec\al,\phi_{\vec\al}},
\eqno(3.57)$$
where
$$\mu_p:\vec\be\mapsto\al_p\be_{_{\ol p}}-\al_{_{\ol p}}\be_p,\;\;\phi_{\vec\al}:\vec\be\mapsto\phi(\vec\al,\vec\be)\eqno(3.58)$$ 
are additive group homomorphisms from $\G$ to $\mbb{F}$. Expressions (3.53)-(3.56) follow from (2.13), (2.14), (3.57) and (3.58).
\psp

Denote
$$\G_3=\{\vec\al\in\rd_\phi\mid\al_{J_{1,3}}=0\},\eqno(3.59)$$
which forms a subgroup of $\G$. 
By Claim 1 and (3.48), we obtain
$$\tau(\G_3)=\G'_3,\eqno(3.60)$$
$$\{\tau(\si_p)\mid p\in \ol{1,\ell_0}\cup I_2\}=\{\si'_p\mid p\in \ol{1,\ell'_0}\cup I'_2\},\eqno(3.61)$$
$$\{\tau(\si_p)\mid p\in \ol{\ell_0+1,\ell_2}\cup I_3\}=\{\si'_p\mid p\in \ol{\ell'_0+1,\ell'_1}\cup I'_3\}.\eqno(3.62)$$
By the above two expressions, 
$$\iota_3=\iota_3'\eqno(3.63)$$
and there exists a bijection $\nu: I_{1,3}\rta I_{1,3}'$ such that
$$\tau(\si_p)=\si'_{\nu(p)} \qquad \for\;\;p\in I_{1,3}.\eqno(3.64)$$

For two subspaces $L_1,L_2$ of ${\cal A}$, we define
$$C_{L_1}(L_2)=\{u\in L_1\,|\,[u,L_2]=0\},\eqno(3.65)$$
$$N_{L_1}(L_2)=\{u\in L_1\,|\,[u,L_2]\subset L_2\}.\eqno(3.66)$$

{\bf Claim 2}. We have
$$C_{\cal A}({\cal A}^{(0)})={\rm Span}\{x^{\vec\al,\vec i}\mid\vec\al\in\rd_\phi,\al_{_{J_{1,3}}}=0,\ \vec i=\vec i_{_{I_6\cup J_7}}\},\eqno(3.67)$$
$$N_{\cal A}({\cal A}^{(0)})={\cal A}^{(0)}+{\rm Span}\{x^{\vec\al,\vec i}\mid
\vec\al\in\rd_\phi,\al_{_{_{J_{1,3}}}}=0,\;|\vec i|=1\mbox{ or }\vec i=\vec i_{_{I_6\cup J_7}}\},\eqno(3.68)$$
where 
$$|\vec i|=\sum_{p\in J}i_p\eqno(3.69)$$
(cf. (2.7), (2.10)).
\psp

Set
$$\es_p=1,\;\;\;\es_{\iota_7+p}=-1\qquad\for\;\;p\in\ol{1,\iota_7}.\eqno(3.70)$$ 
For $p\in J_{1,3}$, we take $0\neq a_{_{[\ol p]}}\in \rd_{\phi}$ by  (2.26) for some $0\neq a\in\mbb{F}$ and have
$$[x^{a_{_{[\ol p]}}},x^{\vec\al,\vec i}]
=\es_{_{\ol p}} a(\al_p x^{\si_p+a_{_{[\ol p]}}+\vec\al,\vec i}+i_p
x^{\si_p+a_{_{[\ol p]}}+\vec\al,\vec i-1_{[p]}})\qquad \for\;(\vec \al,\vec i)\in\G\times {\cal J}.\eqno(3.71)$$
Let
$$u=\sum_{(\vec\al,\vec i)\in\G\times{\cal J}}b_{\vec\al,\vec i}x^{\vec\al,\vec i}\in C_{\cal A}({\cal A}^{(0)})\qquad\mbox{with}\;\;b_{\vec\al,\vec i}\in\mbb{F}.\eqno(3.72)$$
Set
$$k=\mbox{max}\{|\vec i|\mid b_{\vec\al,\vec i}\neq 0\;\mbox{for some}\;\vec\al\in \G\}.\eqno(3.73)$$
By (3.71), we have
$$\al_{_{J_{1,3}}}=0,\;\vec{i}_{_{J_{1,3}}}=0\qquad\mbox{whenever}\;\;b_{\vec\al,\vec i}\neq 0\;\;\mbox{and}\;\;|\vec i|=k.\eqno(3.74)$$
Moreover, by (3.2) and (3.74), we obtain
$$\vec\al\in\rd_{\phi}\qquad\mbox{whenever}\;\;b_{\vec\al,\vec i}\neq 0,\;|\vec i|=k.\eqno(3.75)$$
For $p\in J_4\cup I_{5,6}$, $1_{[p]}\in \G_1$ by (2.15), and we have
$$[x^{1_{[p]}},x^{\vec\be,\vec j}]=\es_pj_{_{\ol p}}x^{\vec\be+1_{[p]},\vec j-1_{[\ol p]}}\eqno(3.76)$$
for $\vec\be\in \rd_{\phi}$ with $\vec\be_{_{J_{1,3}}}=0$ and $\vec j\in{\cal J}$. By (3.2) and (3.74)-(3.76), we get
$$\vec{i}_{_{J_4\cup \ol I_{5,6}}}=0\qquad\mbox{whenever}\;\;b_{\vec\al,\vec i}\neq 0\;\;\mbox{and}\;\;|\vec i|=k.\eqno(3.77)$$
Expressions (3.74), (3.75) and (3.77) imply
$$ \sum_{(\vec\al,\vec i)\in\G\times{\cal J},\;|\vec i|=k}b_{\vec\al,\vec i}x^{\vec\al,\vec i}\in  C_{\cal A}({\cal A}^{(0)}).\eqno(3.78)$$
Hence
$$u-\sum_{(\vec\al,\vec i)\in\G\times{\cal J},\;|\vec i|=k}b_{\vec\al,\vec i}x^{\vec\al,\vec i}=\sum_{(\vec\al,\vec i)\in\G\times{\cal J},\;|\vec i|<k}b_{\vec\al,\vec i}x^{\vec\al,\vec i}\in  C_{\cal A}({\cal A}^{(0)}).\eqno(3.79)$$
By induction on $|\vec i|$, we obtain (3.67). Similarly, we can prove (3.68). 
\psp

Denote
$$
\mbb{F}[\G_3]={\rm Span}\{x^{\vec\al}\mid\vec\al\in\rd_\phi,\al_{_{J_{1,3}}}=0\},
\eqno(3.80)$$
the group algebra of $\G_3$. For $p\in I_{1,3}$, let
$$N_p={\cal A}^{(0)}+\{ u\in N_{\cal A}({\cal A}^{(0)})\mid [x^{\si_p},u]
=0\}.\eqno(3.81)$$
Then $C_{\cal A}({\cal A}^{(0)}),\;N_{\cal A}({\cal A}^{(0)})$ and $N_p$ form $\mbb{F}[\G_3]$-modules.
\psp

{\bf Claim 3}. The quotient  $N_{\cal A}({\cal A}^{(0)})/N_p$  is a free $\mbb{F}[\G_3]$-module of rank 1 with generator
$t_p$ (cf. (2.20)) if $p\in\ol{\ell_0+1,\ell_1}$,  of rank 0 if
$p\in \ol{1,\ell_0}\cup I_2$ and of rank 2 with generators $\{t_p,t_{\ol p}\}$ if $p\in I_3$.
\psp

Set
$${\cal N}= {\rm Span}\{x^{\vec\al,\vec i}\mid
\vec\al\in\rd_\phi,\al_{_{_{J_{1,3}}}}=0,\;|\vec i|=1\mbox{ or }\vec i=\vec i_{_{I_6\cup J_7}}\}.\eqno(3.82)$$
The conclusion follows from (3.68), (3.71) and the fact
$$N_{\cal A}({\cal A}^{(0)})/N_p\cong {\cal N}/\{u\in{\cal N}\mid [x^{\si_p},u]
=0\}.\eqno(3.83)$$ 
\pse

{\bf Claim 4}. We have
$$\nu(\ol{1,\ell_0})=\ol{1,\ell'_0},\;\;\nu(I_i)=I'_i\qquad\for\;\;i=1,2,3.\eqno(3.84)$$
In particular,
$$(\ell_0,\ell_1,\ell_2,\ell_3)=(\ell'_0,\ell'_1,\ell'_2,\ell'_3).\eqno(3.85)$$

By (3.61) and (3.64), we have
$$\nu(I_3)\subset \ol{\ell'_0+1,\ell'_1}\cup I'_3.\eqno(3.86)$$
By (3.48) and (3.64), 
$$\cta(N_{\cal A}({\cal A}^{(0)})/N_p)=N_{{\cal A}'}({{\cal A}'}^{(0)})/N'_{\nu(p)}\qquad\for\;\;p\in I_{1,3}.\eqno(3.87)$$
Moreover, (3.53) implies
$$\cta(\mbb{F}[\G_3])=\mbb{F}[\G'_3].\eqno(3.88)$$
Claim 3 shows that 
$$\nu(I_3)\in I'_3\eqno(3.89)$$
(cf. (2.4)).

Expressions (3.61) and (3.64) imply
$$\nu(\ol{1,\ell_0}\cup I_2)=\ol{1,\ell'_0}\cup I'_2.\eqno(3.90)$$
Suppose that $\nu(p)\in\ol{1,\ell'_0}$ for some $p\in I_2$. Note that
$${\cal N}=\{ u\in N_{\cal A}({\cal A}^{(0)})\mid [u,{\cal A}^{(0)}]\subset\mbb{F}[\G_3]\}\eqno(3.91)$$
by (3.2) (cf. (3.80) and (3.82)). So 
$$\cta({\cal N})={\cal N}'.\eqno(3.92)$$
 Moreover,
$$[x^{\si_p},{\cal N}]=\{0\},\;[{x'}^{\si_{\nu(p)}},{\cal N}']\neq\{0\},\;\;\cta(x^{\si_p})=\chi(\si_p){x'}^{\si_{\nu(p)}}\eqno(3.93)$$
by (3.48).
The above two expressions imply
$$\cta([x^{\si_{[p]}},{\cal N}])\neq [\cta(x^{\si_{[p]}}),\cta({\cal N})],\eqno(3.94)$$
which contradicts that $\cta$ is a Lie algebra isomorphism. Thus
$$\nu(I_2)=I'_2.\eqno(3.95)$$
Since $\nu(I_{1,3})=I'_{1,3}$, we have
 $$\nu(I_1)=I'_1.\eqno(3.96)$$
Furthermore, (3.90) and (3.96) imply
$$\nu(\ol{1,\ell_0})=\ol{1,\ell'_0}.\eqno(3.97)$$
Therefore, Claim 4 holds.
\psp

A linear transformation $T$ on a vector space $V$ is called {\it locally-nilpotent} if
for any $u\in V$, there exist a positive integer $n$ such that $T^n(u)=0$. Set
$${\cal M}=\{ u\in {\cal N}\mid \ad_u|_{{\cal A}^{(0)}}\;\mbox{is locally-nilpotent}\},\eqno(3.98)$$
$${\cal M}_1=\{ u\in {\cal N}\mid \ad_u|_{{\cal A}^{(0)}}\;\mbox{is locally-finite}\}\eqno(3.99)$$
(cf. (1.13), (3.82)).
\psp

{\bf Claim 5}. 
$${\cal M}= {\rm Span}\{x^{\vec\al,\vec i}\mid
\vec\al\in\rd_\phi,\al_{_{_{J_{1,3}}}}=0,\;\vec i=\vec i_{_{I_6\cup J_7}}\}.\eqno(3.100)$$
$${\cal M}_1={\rm Span}\{t_p\,|\,p\in J_4\cup\ol I_{5,6}\}+M.\eqno(3.101)$$
\pse

Note that by (3.2),
$$\ad_u|_{{\cal A}^{(0)}}=0\qquad\for\;\;u\in {\rm Span}\{x^{\vec\al,\vec i}\mid
\vec\al\in\rd_\phi,\al_{_{_{J_{1,3}}}}=0,\;\vec i=\vec i_{_{I_6\cup J_7}}\}.\eqno(3.102)$$
Moreover, for $\al\in\G_3,\;p\in I_1\cup J_{3,4}\cup \ol I_{5,6}$ and $\be\in\G$, we have
$$[x^{\vec\al,1_{[p]}},x^{\vec\be}]=\es_p\be_{_{\ol p}}x^{\si_p+\vec\al+\vec\be}\eqno(3.103)$$
by (3.2). The above two expressions imply (3.100) and (3.101).  
\psp

{\bf Claim 6}. The center of the Lie algebra ${\cal M}$ in (3.100) is
$${\cal M}_0=\sum_{\vec i=\vec i_{I_6}\in{\cal J}}\mbb{F}[\G_3]t^{\vec i}
\eqno(3.104)$$
and
$$\{u\in {\cal M}_1\mid [u,{\cal M}]\subset {\cal M}_0\}=\sum_{p\in J_4\cup \ol I_{5,6}}\mbb{F}t_p+{\cal M}_0.\eqno(3.105)$$
\psp

For $p\in J_7$ and $(\vec\al,\vec i)\in \G\times{\cal J}$, we have
$$
[x^{\vec\al,\vec i},t^2_{_{\ol p}}]=\es_pi_p x^{\vec\al,\vec i-1_{[p]}+1_{[\ol p]}}.\eqno(3.106)$$
By (3.2), we have 
$$[{\cal M}_0,{\cal M}]=\{0\}\eqno(3.107)$$
and 
$$[t_p, x^{\vec\al,\vec i}]=\es_p(\al_{_{\ol p}}x^{\vec\al,\vec i}+i_{_{\ol p}}x^{\vec\al,\vec i-1_{[\ol p]}})\qquad\for\;\;p\in J_4\cup\ol I_{5,6},\;(\vec\al,\vec i)\in \G\times{\cal J}.\eqno(3.108)$$
Expressions (3.104) and (3.105) follow from (3.106)-(3.108).
\psp

By (3.105), we  have
$$\cta(\ol t_{_{J_4\cup\ol I_{5,6}}})
\equiv \ol t'_{_{J'_4\cup\ol I'_{5,6}}}E\ ({\rm mod}\;{\cal M}'_0)\eqno(3.109)$$
(cf. (3.21)) for some invertible matrix
$$E=(e_{p,q})_{(2\ell'_4+\ell'_5+\ell'_6)\times(2\ell_4+\ell_5+\ell_6)}\qquad\mbox{with}\;\;e_{p,q}\in\mbb{F}\eqno(3.110)$$
(cf. (2.12)), where ${\cal M}'_0$ is the corresponding space ${\cal M}_0$ in ${\cal A}'$.
\psp

{\bf Claim 7}. We have
$$e_{p,q}=0\;\mbox{ \ if \ } p\notin\ol I'_5,\;q\in\ol I_5\;\;\mbox{or}\;\;p\in\ol I'_6,\;q\in J_4,\eqno(3.111)$$
which implies
$$(\ell_4,\ell_5,\ell_6)=(\ell'_4,\ell'_5,\ell'_6).\eqno(3.112)$$

Note that 
$$\ad_{t_q}=\es_q\ptl_{_{\ol q}}\qquad\for\;q\in J_4\cup\ol I_{5,6}\eqno(3.113)$$
as operators on ${\cal A}$ by (3.2). Hence $\ad_{t_q}$ is diagonalizable if and only if
$q\in\ol I_5$. But the adjoint operator of the corresponding element $\cta(t_q)$ in
the right-hand side of (3.109) cannot be diagonalizable if $e_{p,q}\ne0$ for some $p\notin\ol I'_5$. So the first case in  (3.111) holds.
By (3.104), 
$$\ad_{t_q}|_{{\cal M}_0}\;\;\mbox{is diagonalizable}\;\;\for\;\;q\in J_4\eqno(3.114)$$
and
$$\cta(t_q)|_{{\cal M}'_0}\;\;\mbox{is not diagonalizable}\;\;\mbox{if}\;e_{p,q}\ne0\;\mbox{for some}\;p\in \ol I'_6.\eqno(3.115)$$
The above two expressions imply  the second case in (3.111).
\psp

We have
$$(\ell_0,\vec\ell)=(\ell'_0,\vec\ell').\eqno(3.116)$$ 
by Claims 4, 7 and the fact $\cta({\cal M})={\cal M}'$, where ${\cal M}'$ is the subalgebra of ${\cal A}'$ corresponding to the subalgebra ${\cal M}$ of ${\cal A}$ in (3.98) (also cf. (3.100)). Moreover, (3.111) implies that the matrix
$$E=\left(\begin{array}{ccc} A_1'&0&C'\\ B'&A'_2&D'\\ 0&0&A'_3\end{array}\right),\eqno(3.117)$$
where 
$$  A'_1\in GL_{2\ell_4}(\mbb{F}),\;\;A'_2\in GL_{\ell_5}(\mbb{F}),\;\;A'_3\in GL_{\ell_6}(\mbb{F}),\eqno(3.118)$$
$$B'\in M_{\ell_5\times 2\ell_4}(\mbb{F}),\;\;C'\in M_{2\ell_4\times \ell_6}(\mbb{F}),\;\;D'\in M_{\ell_5\times \ell_6}(\mbb{F}).\eqno(3.119)$$
In particular,
$$\cta(\ol t_{_{J_4}})\equiv \ol t'_{_{J_4}}A'_1+\ol t'_{_{\ol I_5}}B'\;\; ({\rm mod}\;{\cal M}'_0).\eqno(3.120)$$
Set
$$\widetilde{\cal M}'=\sum_{0\neq \vec\al\in\G_3,\;\vec i=\vec i_{I_6}\in{\cal J}}
\mbb{F}x^{\vec\al,\vec i}.\eqno(3.121)$$
By (3.2), (3.59) and (3.104), we have
$$[\ol t'_{_{J_4\cup\ol I_5}},{\cal M}'_0]\subset \widetilde{\cal M}'.\eqno(3.122)$$
Thus we have
\begin{eqnarray*}\hspace{1cm}S_{\ell_4}&=&\cta([\ol t^T_{_{J_4}},\ol t_{_{J_4}}])\\ &=&[\cta(\ol t_{_{J_4}})^T,\cta(\ol t_{_{J_4}})]\\ &\equiv& [(\ol t'_{_{J_4}}A'_1+\ol t'_{_{\ol I_5}}B')^T,\ol t'_{_{J_4}}A'_1+\ol t'_{_{\ol I_5}}B']\;\; ({\rm mod}\;\widetilde{\cal M}')\\ &\equiv& {A_1'}^TS_{\ell_4}A_1'\;\; ({\rm mod}\;\widetilde{\cal M}')\hspace{8.4cm}(3.123)\end{eqnarray*}
(cf. (3.6)). The connection of $E$ with the $f$ in (3.8) is
$$E^{-1}=\left(\begin{array}{ccc} A_1&0&C\\ B&A_2&D\\ 0&0&A_3\end{array}\right).\eqno(3.124)$$

Denote by $\pi'_0$ and $\pi'_1$ the projections from $\G'$ to $\G'_0$ and $\G'_1$, respectively. Set
$$\tau_1=\pi'_0\tau|_{\G_1}:\G_1\rta\G'_0,\;\;\tau_2=\pi'_1\tau|_{\G_1}:\G_1\rta\G'_1.\eqno(3.125)$$
Then $\tau_1$ and $\tau_2$ are group homomorphisms. We shall determine $\tau_2$. Let $\vec\al\in\G$ be an arbitrary element. For  $p\in I_{1,3}$, we pick 
$$\vec\g=a_{_{[\ol p]}}\in(\mbb{F}1_{[\ol p]}\cap\rd_\phi)\setminus\{0\}\eqno(3.126)$$
by (2.26). By (3.2), (3.48) and (3.64)
\begin{eqnarray*}\hspace{1.5cm}& &\chi(\vec\al)\chi(\vec\g)(\phi'(\tau(\vec\al),\tau(\vec\g)){x'}^{\tau(\vec\al)+\tau(\vec\g)}+
\sum_{q\in I_{1,3}}((\tau_2(\al))_q(\tau_2(\vec\g))_{_{\ol q}}-\\& &-(\tau_2(\al))_{_{\ol q}}(\tau_2(\vec\g))_q) {x'}^{\si'_q+\tau(\vec\al)+\tau(\vec\g)}) \\&=&\chi(\vec\al)\chi(\vec\g) [{x'}^{\tau(\vec\al)},{x'}^{\tau(\vec\g)}]\\ &=&[\cta(x^{\vec\al}),\cta(x^{\vec\g})]\\ &=&\cta([x^{\vec\al},x^{\vec\g}])\\ &=&\al_p
a\cta(x^{\si_p+\vec\al+\vec\g})\\&=&\al_p
a\chi(\si_p+\vec\al+\vec\g) {x'}^{\si'_{\nu(p)}+\tau(\vec\al)+\tau(\vec\g)}.\hspace{6.2cm}(3.127)\end{eqnarray*}
Comparing the coefficients of ${x'}^{\si'_{\nu(p)}+\tau(\vec\al)+\tau(\vec\g)}$, we obtain
$$(\tau_2(\al))_{\nu(p)}(\tau_2(\vec\g))_{_{\ol {\nu(p)}}}-(\tau_2(\al))_{_{\ol {\nu(p)}}}(\tau_2(\vec\g))_{\nu(p)}=a\chi(\si_p)\al_p.\eqno(3.128)$$

First suppose $p\in\ol{1,\ell_0}$. By (3.1), (3.48) and (3.64),
\begin{eqnarray*}\hspace{3cm}& &-\al_{_{\ol p}}\chi(\vec\al){x'}^{\tau(\vec\al)}\\&=&
-\al_{_{\ol p}}\cta(x^{\vec\al})\\ &=&\cta([x^{-\si_p},x^{\vec\al}])\\ &=&[\cta(x^{-\si_p}),\cta(x^{\vec\al})]\\ &=&\chi(-\si_p)\chi(\vec\al)[{x'}^{-\si'_{\nu(p)}},{x'}^{\tau(\vec\al)}]\\ &=&-\chi(-\si_p)\chi(\vec\al)(\tau_2(\al))_{_{\ol{\nu(p)}}}{x'}^{\tau(\vec\al)}.\hspace{5.4cm}(3.129)\end{eqnarray*}
This shows
$$(\tau_2(\al))_{_{\ol{\nu(p)}}}=\chi(\si_p)\al_{_{\ol p}}\qquad\for\;\;\al\in \G_1,\;p\in\ol{1,\ell_0}.\eqno(3.130)$$
By (3.126) and the above expression,
$$(\tau_2(\vec\g))_{_{\ol \nu(p)}}=a\chi(\si_p).\eqno(3.131)$$
Substitute (3.130) and (3.131) into (3.128):
$$a\chi(\si_p)(\tau_2(\al))_{\nu(p)}-\chi(\si_p)\al_{_{\ol p}}(\tau_2(\vec\g))_{\nu(p)}=a\chi(\si_p)\al_p ,\eqno(3.132)$$
which implies
$$(\tau_2(\al))_{\nu(p)}=\al_p+a^{-1}(\tau_2(\vec\g))_{\nu(p)}\al_{_{\ol p}}\qquad\for\;\;\al\in \G_1,\;p\in\ol{1,\ell_0}.\eqno(3.133)$$
Expressions (3.130) and (3.132) relate to the first equation in (3.23) with
$$a_p=a^{-1}(\tau_2(\vec\g))_{\nu(p)},\;\;b_p=\chi(\si_p)\qquad\for\;\;p\in\ol{1,\ell_0}.\eqno(3.134)$$

Next suppose $p\in \ol{\ell_0+1,\ell_1}\cup I_{2,3}$. By (3.1), (3.48) and (3.64),
\begin{eqnarray*}\hspace{2cm}& &(\al_p-\al_{_{\ol p}})\chi(\vec\al){x'}^{\tau(\vec\al)}\\&=&(\al_p-\al_{_{\ol p}})\cta(x^{\vec\al})\\ &=&\cta([x^{-\si_p},x^{\vec\al}])\\ &=&[\cta(x^{-\si_p}),\cta(x^{\vec\al})]\\ &=&\chi(-\si_p)\chi(\vec\al)[{x'}^{-\si'_{\nu(p)}},{x'}^{\tau(\vec\al)}]\\ &=&\chi(-\si_p)\chi(\vec\al)((\tau_2(\al))_{\nu(p)}-(\tau_2(\al))_{_{\ol{\nu(p)}}})
{x'}^{\tau(\vec\al)}.\hspace{4cm}(3.135)\end{eqnarray*}
Thus we have
$$\al_p-\al_{_{\ol p}}=\chi(-\si_p)((\tau_2(\al))_{\nu(p)}-(\tau_2(\al))_{_{\ol{\nu(p)}}}),\eqno(3.136)$$
which implies
$$(\tau_2(\al))_{\nu(p)}-(\tau_2(\al))_{_{\ol{\nu(p)}}}=\chi(\si_p)(\al_p-\al_{_{\ol p}}).\eqno(3.137)$$
In particular,
$$(\tau_2(\vec\g))_{\nu(p)}-(\tau_2(\vec\g))_{_{\ol{\nu(p)}}}=-a\chi(\si_p).\eqno(3.138)$$
Solving (3.128) by (3.137) and (3.138), we obtain 
$$((\tau_2(\al))_{\nu(p)},(\tau_2(\al))_{_{\ol{\nu(p)}}})=(\al_p,\al_{_{\ol p}})\left(\begin{array}{cc}a_p+b_p&a_p\\ 1-a_p-b_p&1-a_p\end{array}\right),
\eqno(3.139)$$
where 
$$a_p=1-a^{-1}(\tau_2(\vec\g))_{_{\ol{\nu(p)}}},\;\;b_p=\chi(\si_p)\ne 0.\eqno(3.140)$$
This gives the third equation in (3.23) by (3.13).

 Assume $p\in\ol{\ell_0+1,\ell_1}$. Set
$$\bar\G_p'=(\mbb{F}1_{[\nu(p)]}+\mbb{F}1_{[\ol{\nu(p)}]})\bigcap \rd_{\phi'},\;\;\bar\G_p=\tau^{-1}(\bar\G_p').\eqno(3.141)$$
By (2.26), $\bar\G_p'$ contains a basis of $\mbb{F}1_{[\nu(p)]}+\mbb{F}1_{[\ol{\nu(p)}]}$. Note that
$$\left|\begin{array}{cc}a_p+b_p&a_p\\ 1-a_p-b_p&1-a_p\end{array}\right|=b_p\neq 0,\eqno(3.142)$$
which implies that the matrix
$$\left(\begin{array}{cc}a_p+b_p&a_p\\ 1-a_p-b_p&1-a_p\end{array}\right)\;\;\mbox{is invertible}.\eqno(3.143)$$
Thus
\begin{eqnarray*}\hspace{1.5cm}& &\{(\al_p,\al_{_{\ol p}})\mid \vec\al=(\al_0,\al_1,\al_{_{\ol 1}},...,\al_{\iota_7},\al_{_{\ol{\iota_7}}})\in \bar\G_p\}\\ & &\mbox{is an additive subgroup of}\;\;\mbb{F}^{\:2}\;\;\mbox{containing a basis of}\;\;\mbb{F}^{\:2}\hspace{2.3cm}(3.144)\end{eqnarray*}
by (3.139) and (3.141). Let $\vec\al\in\bar\G_p$. Note that
\begin{eqnarray*}[{x'}^{\vec\be',\vec j},{x'}^{\tau(\vec\al)}]&=&(\be'_{\nu(p)}(\tau(\vec\al))_{_{\ol{\nu(p)}}}-\be'_{_{\ol{\nu(p)}}}(\tau(\vec\al))_{\nu(p)}){x'}^{\si'_{\nu(p)}+\vec \be'+\tau(\vec\al),\vec{j}}\\ & &+j_{\nu(p)}(\tau(\vec\al))_{_{\ol{\nu(p)}}}{x'}^{\si'_{\nu(p)}+\vec \be'+\tau(\vec\al),\vec{j}-1_{[\nu(p)]}}\hspace{5.2cm}(3.145)\end{eqnarray*}
by (3.2) and the fact $\tau(\vec\al)\in\bar\G'_p$  for $(\vec\be',\vec j)\in \G'\times {\cal J}'$. Furthermore,
\begin{eqnarray*}\hspace{3cm}\al_{_{\ol p}}\chi(\vec\al+\si_p){x'}^{\si'_{\nu(p)}+\tau(\vec{\al})}&=&\cta(\al_{_{\ol p}}x^{\si_p+\vec{\al}})\\ &=&\cta([t_p,x^{\vec\al}])\\ &=& [\cta(t_p),\cta(x^{\vec\al})]\\ &=&\chi(\vec\al)[\cta(t_p),{x'}^{\tau(\vec{\al})}]\hspace{3.6cm}(3.146)\end{eqnarray*}
by  (3.48). We write
$$\cta(t_p)=bt'_{\nu(p)}+\sum_{(0,1_{\nu(p)})\neq (\vec\be',\vec j)\in \G'\times{\cal J}'}b_{\vec\be',\vec j}{x'}^{\vec\be',\vec j}.\eqno(3.147)$$
Expressions (3.145)-(3.147) show
$$b(\tau(\vec\al))_{_{\ol{\nu(p)}}}=\chi(\si_p)\al_{_{\ol p}}=b_p\al_{_{\ol p}}.\eqno(3.148)$$
Furthermore, (3.139) and (3.148) imply
$$b(a_p\al_p+(1-a_p)\al_{\ol p})=b_p\al_{\ol p},
\eqno(3.149)$$
which implies
$$a_p=0\eqno(3.150)$$
by (3.144). This gives  the second equation in (3.23) by (3.13).

Note that ${\cal M}_0$ in (3.104) commutes with ${\cal A}^{(0)}$ in (3.49). By (3.2) and (3.109), we have
\begin{eqnarray*}\hspace{1.6cm}\chi(\vec\al)\al_{_{J_4\cup I_{5,6}}}{x'}^{\tau(\vec\al)}&=& \cta(\al_{_{J_4\cup I_{5,6}}}x^{\vec\al})\\ &=&\cta([\ol t_{_{J_4\cup\ol I_{5,6}}},x^{\vec\al}])\\&=& [\cta(\ol t_{_{J_4\cup\ol I_{5,6}}}),\cta(x^{\vec\al})]\\ &=&\chi(\vec\al)[\ol t'_{_{J_4\cup\ol I_{5,6}}},{x'}^{\tau(\vec\al)}]E
\\ &=&\chi(\vec\al)(\tau(\vec\al))_{_{J_4\cup I_{5,6}}}E{x'}^{\tau(\vec\al)},\hspace{4.7cm}(3.151)\end{eqnarray*}
Thus 
$$\al_{_{J_4\cup I_{5,6}}}=(\tau(\vec\al))_{_{J_4\cup I_{5,6}}}E.\eqno(3.152)$$
This gives the action of $f$ in (3.13), which is determined by (3.8) with the corresponding data in (3.124). This shows that $\tau_2\in G$ in (3.14).
In particular,  $\tau_2$ must be an isomorphism. If
 $\vec\al\in\G_0$, then $\al_p=0$ and so $\tau(\vec\al)_p=0$ for all
$p\in J$. Hence 
$$\tau_0=\tau|_{\G_0}:\G_0\rta\G'_0\eqno(3.153)$$
 is an isomorphism. Applying $\cta$ to (3.2), we obtain
that 
$$\phi'(\tau(\vec\al),\tau(\vec\be))=\phi(\vec\al,\vec\be)\qquad\for\;\;\vec\al,\vec\be\in\G\eqno(3.154)$$
by (3.13). Moreover, (3.15) holds by (3.125) and (3.153). Equations (3.16) and (3.17) follow from (3.15) and (3.154). This completes the proof of the Theorem.$\qquad\Box$

\vspace{0.7cm}

\noindent{\Large \bf References}

\hspace{0.3cm}

\begin{description}

\item[{[K1]}] V. G. Kac, A description of filtered Lie algebras whose associated graded Lie algebras are of Cartan types, {\it Math. of USSR-Izvestijia} {\bf 8} (1974), 801-835.

\item[{[K2]}] V. G. Kac, Lie superalgebras, {\it Adv. Math.} {\bf 26} (1977), 8-96.

\item[{[K3]}] V. G. Kac, Classification of infinite-dimensional simple linearly compact Lie superalgebras, {\it Adv. Math.} {\bf 139} (1998), 1-55.

\item[{[O]}] J. Marshall Osborn, New simple infinite-dimensional Lie algebras of characteristic 0, {\it J. Algebra} {\bf 185} (1996), 820-835.

\item[{[OZ]}]   J. Marshall Osborn and K. Zhao, Generalized Poisson brackets and Lie algebras for type $H$ in characteristic 0, {\it Math. Z.} {\bf 230} (1999), 107-143. 

\item[{[SXZ]}] Y. Su, X. Xu and H. Zhang, Derivation-simple algebras and the structures of  Lie algebras of Witt type, {\it J. Algebra} {\bf 233} (2000), 642-662. 

\item[{[X]}] X. Xu, New generalized simple Lie algebras of Cartan type over a field with characteristic 0, {\it J. Algebra} {\bf 224} (2000), 23-58.

\end{description}
\end{document}